\documentclass[11pt,a4paper]{article}

\usepackage[T1]{fontenc}
\usepackage[utf8]{inputenc}
\usepackage{lmodern}
\usepackage{amsmath,amssymb,amsthm,mathtools}
\usepackage{booktabs}
\usepackage[margin=2.5cm]{geometry}
\usepackage[colorlinks=true,linkcolor=blue,citecolor=blue,urlcolor=blue,
            bookmarksnumbered=true,
            pdftitle={On the fourth power mean of one kind three-term character sums},
            pdfauthor={Hu Jiayuan and Tang Jun}]{hyperref}

\DeclareMathOperator{\ord}{ord}
\newcommand{\ep}{e_{p}}
\newcommand{\Fp}{\mathbb{F}_{p}}
\newcommand{\Fps}{\mathbb{F}_{p}^{*}}
\newcommand{\chibar}{\overline{\chi}}
\newcommand{\ind}{\mathbf{1}}

\theoremstyle{plain}
\newtheorem{theorem}{Theorem}
\newtheorem{lemma}{Lemma}
\newtheorem{corollary}{Corollary}
\theoremstyle{remark}
\newtheorem{remark}{Remark}

\allowdisplaybreaks

\begin{document}

\title{\bfseries On the fourth power mean of\\[2pt]
one kind three-term character sums%
  \thanks{This work is supported by the National Natural Science Foundation of China
  (Grant No.~12161031) and the Natural Science Foundation of Inner Mongolia (CN)
  (Grant No.~2026MS0472).}}
\author{Hu Jiayuan$^{1}$ and Tang Jun$^{2,\dagger}$\\[4pt]
  {\footnotesize $^{1}$School of Science, Inner Mongolia University of Science and Technology,}\\[-3pt]
  {\footnotesize Baotou, Inner Mongolia, P. R. China}\\[1pt]
  {\footnotesize $^{2}$School of Energy and Environment, Inner Mongolia University of Science}\\[-3pt]
  {\footnotesize and Technology, Baotou, Inner Mongolia, P. R. China}}
\date{}
\maketitle

\begingroup
\renewcommand{\thefootnote}{\fnsymbol{footnote}}
\footnotetext[2]{Corresponding author. E-mail: \texttt{tj2001766@sina.com}.}
\endgroup

\noindent\textbf{Abstract}

\noindent The main purpose of this paper is to give an exact formula for the fourth power
mean of one kind three-term character sums mod $p$, a problem left open by Zhang and Hu
(2017). For any non-principal character $\chi$ mod $p$ of order $d$, we prove that the
mean depends only on $p$ and $d$, and determine it in the four cases $d=2,3,4$ and
$d\notin\{2,3,4\}$. The proof uses Gauss sum inversion, orthogonality and a Vandermonde
identity, and we also give the asymptotic behaviour and a closed formula for the
principal character.



%=====================================================================
\section{Introduction}
%=====================================================================

Let $p\ge 3$ be a prime, $\ep(t)=e^{2\pi it/p}$, and let $\chi$ denote a Dirichlet
character mod $p$ with the convention $\chi(0)=0$. This convention also applies to the
principal character $\chi_{0}$ and to all powers of characters. For any integers $m$ and
$n$, the three-term character sum $C(m,n;\chi;p)$ and its fourth power mean
$f_{\chi}(p)$ are defined as follows:
\begin{equation}\label{eq:def}
C(m,n;\chi;p)=\sum_{a=1}^{p-1}\chi\bigl(ma^{3}+na^{2}+a\bigr),
\qquad
f_{\chi}(p)=\sum_{m=0}^{p-1}\sum_{n=0}^{p-1}\bigl|C(m,n;\chi;p)\bigr|^{4}.
\end{equation}
If one wants the mean value with respect to $(m,n)$, it is enough to divide by $p^{2}$.
We keep the unnormalized sum \eqref{eq:def}, which is the convention used in [6]--[10].

Note that
\[
ma^{3}+na^{2}+a=a\bigl(ma^{2}+na+1\bigr),
\]
so the polynomial acted on by $\chi$ has constant term $0$ and the coefficient of the
linear term equal to $1$. This is the essential difference between the three-term
character sum \eqref{eq:def} and the classical two-term exponential sum
$\sum_{a}\ep(ma^{k}+na)$.

Upper bounds for character sums and exponential sums are basic tools in analytic number
theory. Weil [3] proved that if $\chi$ is a character of order $d>1$ mod $p$ and
$F\in\Fp[x]$ is a polynomial of positive degree that is not a constant multiple of a
$d$-th power, then
\[
\Bigl|\sum_{x\in\Fp}\chi\bigl(F(x)\bigr)\Bigr|\le(\deg F-1)\sqrt{p}.
\]
Burgess [4] studied the distribution of quadratic residues and non-residues, while
Cochrane and Zheng [5] obtained bounds for related exponential sums. However, the
distribution of a single sum is very irregular, while its higher power mean often
satisfies an exact identity. This phenomenon is the main motivation of a series of works
of Zhang Wenpeng and his coauthors; see [6]--[9].

Closely related to the present paper is the hybrid power mean of Zhang Wenpeng and Hu
Jiayuan [6]. Let $p$ be an odd prime with $(3,p-1)=1$. Then for any non-principal
character $\chi$ mod $p$ they proved that
\begin{equation}\label{eq:hybrid}
\begin{aligned}
&\sum_{m=0}^{p-1}\sum_{n=0}^{p-1}
\Bigl|\sum_{a=1}^{p-1}\chi\bigl(ma^{3}+na^{2}+a\bigr)\Bigr|^{2}
\Bigl|\sum_{b=1}^{p-1}\ep\bigl(mb^{3}+nb^{2}\bigr)\Bigr|^{2}\\
&\qquad=\begin{cases}
p(p-1)(p^{2}-2p+2), & \chi\ \text{non-real odd},\\[1pt]
p(p^{3}-5p^{2}+10p-2), & \chi\ \text{non-real even},\\[1pt]
p^{2}(p-1)(p-2), & \chi\ \text{Legendre},\ p\equiv3\pmod{4},\\[1pt]
p^{2}(p^{2}-5p+8), & \chi\ \text{Legendre},\ p\equiv1\pmod{4}.
\end{cases}
\end{aligned}
\end{equation}
In \eqref{eq:hybrid} the character sum occurs only to the second power. In the concluding
remarks of [6] the authors explicitly noted that they could not give an exact
computational formula for the pure fourth power mean
$\sum_{m,n}\bigl|\sum_{a}\chi(ma^{3}+na^{2}+a)\bigr|^{4}$, and listed it as a problem
for further research. The present paper solves this problem by evaluating $f_{\chi}(p)$
exactly, so the two results are complementary. Exact
formulae for the fourth power mean of three-term exponential sums were obtained by Liu
Huaning and Li Wanmei [7] with the help of trigonometric identities and the properties of
reduced residue systems, by Du Xiancun and Han Di [9] and by Liu and Chen [10] for other
pairs of exponents; Bai and Hu [11] discussed the error term of the fourth power mean of
generalized cubic Gauss sums. The method used below combines Gauss sum expansion,
orthogonality, and the evaluation of a Jacobi sum. An additional coefficient is introduced
so that the orthogonality conditions form a Vandermonde system.

A direct two-parameter reduction is obtained by writing $ma^{3}+na^{2}+a=a\,h(a)$
with $h(t)=mt^{2}+nt+1$, expanding $|C|^{4}$ and applying the Gauss sum inversion
formula to the four factors. This gives
\begin{equation}\label{eq:thm1}
f_{\chi}(p)=\sum_{r_{1},r_{2},r_{3},r_{4}\in\Fps}
\chibar(r_{1}r_{2})\chi(r_{3}r_{4})\,\ep(r_{1}+r_{2}+r_{3}+r_{4})\,N(r_{1},r_{2},r_{3},r_{4}),
\end{equation}
where
\begin{equation}\label{eq:N}
N(r)=\sum_{\substack{a,b,c,d\in\Fps\\ r_{1}a+r_{2}b+r_{3}c+r_{4}d=0\\
r_{1}a^{2}+r_{2}b^{2}+r_{3}c^{2}+r_{4}d^{2}=0}}
\chi(a)\chi(b)\chibar(c)\chibar(d).
\end{equation}
All equations involving elements of $\Fp$ are understood in that field, and a derivation
of \eqref{eq:thm1}--\eqref{eq:N} is given at the beginning of Section~3. The phase in
\eqref{eq:thm1} does not involve $a,b,c,d$, but the variables in \eqref{eq:N} remain
coupled by the two congruences; this is precisely why the exact evaluation below
introduces a third coefficient.

The main result of this paper is the following theorem, which gives four exact formulae
for $f_{\chi}(p)$ according to the order of $\chi$. Corollary~1 below gives the asymptotic
behaviour and exhibits the jump of the leading coefficient.

\begin{theorem}
Let $p$ be an odd prime, $\chi$ a non-principal character mod $p$ and
$d=\ord(\chi)$. Then $f_{\chi}(p)$ depends only on $d$ and $p$, and
\begin{equation}\label{eq:main}
f_{\chi}(p)=
\begin{cases}
2p(p-1)(p^{2}-p-3)=2p^{4}-4p^{3}-4p^{2}+6p, & d=2,\\[1pt]
p(p-1)(p^{2}-7)=p^{4}-p^{3}-7p^{2}+7p, & d=3,\\[1pt]
3p^{2}(p-1)(p-3)=3p^{4}-12p^{3}+9p^{2}, & d=4,\\[1pt]
2p^{2}(p-1)(p-4)=2p^{4}-10p^{3}+8p^{2}, & d\notin\{2,3,4\}.
\end{cases}
\end{equation}
The case $d=2$ means that $\chi$ is the Legendre symbol mod $p$; the case $d=3$ occurs
only when $3\mid(p-1)$ and the case $d=4$ only when $4\mid(p-1)$, so the corresponding
branches are absent for $p\equiv2\pmod{3}$ or $p\equiv3\pmod{4}$.
\end{theorem}

\begin{corollary}
With the notation of Theorem 1, assume in addition that $d\notin\{2,3,4\}$. Then,
uniformly over such characters,
\begin{equation}\label{eq:cor}
f_{\chi}(p)=2p^{4}-10p^{3}+O(p^{2}),\qquad p\to\infty.
\end{equation}
Equivalently, $p^{-2}f_{\chi}(p)=2p^{2}-10p+O(1)$. The leading coefficients of the three
exceptional branches are $2$, $1$, $3$ for $d=2,3,4$, respectively: an order three
character halves the main term, while an order four character enlarges it by a factor of
$1.5$.
\end{corollary}

\begin{remark}
The four identities of Theorem 1 agree with numerical computations for every odd prime
$p\le97$ and every non-principal character mod $p$; see Section 4. Their proof in Section
3 does not depend on these computations. Instead of assuming cancellation in the
constrained sum \eqref{eq:N}, we evaluate an enlarged moment in which the linear coefficient also
varies, and then subtract its zero-coefficient part.
\end{remark}

%=====================================================================
\section{Several Lemmas}
%=====================================================================

In this section we give the lemmas needed below. We use the standard properties of Gauss
sums, including $|\tau(\psi)|=\sqrt{p}$ and $\tau(\psi)\tau(\overline{\psi})=\psi(-1)p$
for non-principal $\psi$; see [1,2].

\begin{lemma}[Gauss sum inversion]\label{lem:gauss}
Let $\chi$ be a non-principal character mod $p$ and
$\tau(\chi)=\sum_{a=1}^{p-1}\chi(a)\ep(a)$. Then for any integer $u$, including
$u\equiv0\pmod{p}$,
\begin{equation}\label{eq:inv}
\chi(u)=\frac{1}{\tau(\chibar)}\sum_{r=1}^{p-1}\chibar(r)\ep(ru),
\qquad
\chibar(u)=\frac{1}{\tau(\chi)}\sum_{r=1}^{p-1}\chi(r)\ep(ru).
\end{equation}
\end{lemma}

\begin{proof}
We prove the first identity; the second follows on replacing $\chi$ by $\chibar$. If
$u\not=0$, put $t=ru$. Then
\[
\sum_{r=1}^{p-1}\chibar(r)\ep(ru)=\sum_{t=1}^{p-1}\chibar\bigl(tu^{-1}\bigr)\ep(t)
=\chi(u)\tau(\chibar).
\]
If $u=0$, the sum on the right of \eqref{eq:inv} vanishes because $\chi$ is non-principal, in
agreement with $\chi(0)=0$. This proves Lemma \ref{lem:gauss}.
\end{proof}

\begin{lemma}[Orthogonality]\label{lem:orth}
For $X,Y,R\in\Fp$,
\begin{equation}\label{eq:orth}
\sum_{m=0}^{p-1}\sum_{n=0}^{p-1}\ep(mX+nY+R)=p^{2}\ep(R)\ind_{X=0}\ind_{Y=0}.
\end{equation}
\end{lemma}

\begin{proof}
Each geometric sum is $p$ when its coefficient is zero and is zero otherwise. Their
product gives \eqref{eq:orth}. The same argument gives
\[
\sum_{m,n,\ell\in\Fp}\ep(mX+nY+\ell Z)=p^{3}\ind_{X=0}\ind_{Y=0}\ind_{Z=0}.
\qedhere
\]
\end{proof}

\begin{lemma}[Weil bound]\label{lem:weil}
Let $\chi$ have order $d>1$ mod $p$, and let $F\in\Fp[x]$ have positive degree. Assume
that $F$ is not of the form $cG^{d}$ with $c\in\Fps$ and $G\in\Fp[x]$. Then
\[
\Bigl|\sum_{x\in\Fp}\chi\bigl(F(x)\bigr)\Bigr|\le(\deg F-1)\sqrt{p}.
\]
\end{lemma}

\begin{proof}
This is the classical multiplicative character form of the Weil bound; see [2, Ch.~11]
and [3].
\end{proof}

\begin{lemma}[A Jacobi sum]\label{lem:jacobi}
Put
\begin{equation}\label{eq:JA}
J=J(\chi,\chi^{2})=\sum_{t\in\Fp}\chi(t)\chi^{2}(1-t),\qquad A=|J|^{2}.
\end{equation}
If $\chi$ is non-principal of order $d$, then
\begin{equation}\label{eq:Aval}
A=\begin{cases} 1, & d=2\ \text{or}\ 3,\\[1pt] p, & d\notin\{2,3\}.\end{cases}
\end{equation}
\end{lemma}

\begin{proof}
If $d=2$, then $\chi^{2}=\chi_{0}$, with $\chi_{0}(0)=0$, and
\[
J=\sum_{t\not=1}\chi(t)=-1.
\]
If $d=3$, then $\chi^{2}=\chibar$. The substitution $u=t/(1-t)$ maps
$\Fp\setminus\{1\}$ bijectively onto $\Fp\setminus\{-1\}$, whence
\[
J=\sum_{t\not=1}\chi\Bigl(\frac{t}{1-t}\Bigr)=-\chi(-1).
\]
In either case $A=1$. Otherwise $\chi$, $\chi^{2}$ and $\chi^{3}$ are all non-principal.
Expanding $\tau(\chi)\tau(\chi^{2})$ and grouping by $s=u+v$, the part $s=0$ is
$\chi^{2}(-1)\sum_{u\not=0}\chi^{3}(u)=0$. For $s\not=0$, put $u=st$ and $v=s(1-t)$. It
follows that
\[
\tau(\chi)\tau(\chi^{2})=J\tau(\chi^{3}).
\]
Taking squared absolute values gives $A=p$. This proves Lemma \ref{lem:jacobi}.
\end{proof}

\begin{lemma}[Evaluation maps and a Vandermonde kernel]\label{lem:vander}
Put $F_{m,n,\ell}(t)=mt^{3}+nt^{2}+\ell t$. For $k\le3$ distinct nonzero points
$x_{1},\dots,x_{k}$, the map
\[
(m,n,\ell)\longmapsto\bigl(F_{m,n,\ell}(x_{1}),\dots,F_{m,n,\ell}(x_{k})\bigr)
\]
is onto $\Fp^{k}$, and every fibre has $p^{3-k}$ elements. For four distinct nonzero
points $x_{1},x_{2},x_{3},x_{4}$, put $P(X)=\prod_{i=1}^{4}(X-x_{i})$. All solutions of
\[
\sum_{i=1}^{4}r_{i}x_{i}=\sum_{i=1}^{4}r_{i}x_{i}^{2}=\sum_{i=1}^{4}r_{i}x_{i}^{3}=0
\]
are
\begin{equation}\label{eq:kernel}
r_{i}=\frac{\lambda}{x_{i}P'(x_{i})},\qquad 1\le i\le4,\quad \lambda\in\Fp.
\end{equation}
In particular, all $r_{i}$ are nonzero exactly when $\lambda\not=0$.
\end{lemma}

\begin{proof}
The evaluation matrix with columns ordered as $\ell,n,m$ has rows $(x_{i},x_{i}^{2},x_{i}^{3})$.
Its first $k$ columns have determinant
$\bigl(\prod_{i=1}^{k}x_{i}\bigr)\prod_{i<j}(x_{j}-x_{i})\not=0$. This proves the first
assertion and shows that the three by four matrix in the second assertion has rank three.
Lagrange interpolation, applied to $1$, $X$, $X^{2}$, gives
\[
\sum_{i=1}^{4}\frac{x_{i}^{j}}{P'(x_{i})}=0,\qquad j=0,1,2,
\]
by comparison of coefficients of $X^{3}$. Thus \eqref{eq:kernel} lies in the kernel; since the kernel is
one-dimensional, it gives all solutions.
\end{proof}

%=====================================================================
\section{Proofs of the theorems}
%=====================================================================

Let $p$ be an odd prime and $\chi$ a non-principal character mod $p$. We first retain the
two-parameter reduction, and then evaluate the fourth moment by allowing the linear
coefficient to vary.

\subsection*{Derivation of the two-parameter reduction}

Put $h(t)=mt^{2}+nt+1$. Since $\chi$ is completely multiplicative, including at zero,
\begin{equation}\label{eq:chi}
\chi\bigl(ma^{3}+na^{2}+a\bigr)=\chi(a)\chi\bigl(h(a)\bigr).
\end{equation}
Writing $|C|^{4}=C\,C\,\overline{C}\,\overline{C}$ and summing over $a,b,c,d$, we obtain
\begin{equation}\label{eq:expand}
\begin{aligned}
f_{\chi}(p)&=\sum_{m,n\in\Fp}\Bigl|\sum_{a\in\Fps}\chi(a)\chi\bigl(h(a)\bigr)\Bigr|^{4}\\
&=\sum_{a,b,c,d\in\Fps}\chi(a)\chi(b)\chibar(c)\chibar(d)
\sum_{m,n\in\Fp}\chi\bigl(h(a)h(b)\bigr)\chibar\bigl(h(c)h(d)\bigr).
\end{aligned}
\end{equation}
Applying \eqref{eq:inv} to the four factors and using
$\tau(\chibar)^{2}\tau(\chi)^{2}=(\chi(-1)p)^{2}=p^{2}$, we get
\begin{equation}\label{eq:inv4}
\chi\bigl(h(a)h(b)\bigr)\chibar\bigl(h(c)h(d)\bigr)
=\frac{1}{p^{2}}\sum_{r_{1},r_{2},r_{3},r_{4}\in\Fps}
\chibar(r_{1}r_{2})\chi(r_{3}r_{4})\,\ep\Bigl(\sum_{i=1}^{4}r_{i}h(t_{i})\Bigr),
\end{equation}
where $(t_{1},t_{2},t_{3},t_{4})=(a,b,c,d)$. Since $h(t)=mt^{2}+nt+1$,
\begin{equation}\label{eq:XNY}
\sum_{i=1}^{4}r_{i}h(t_{i})=mX+nY+R,\qquad
X=\sum_{i=1}^{4}r_{i}t_{i}^{2},\quad
Y=\sum_{i=1}^{4}r_{i}t_{i},\quad
R=\sum_{i=1}^{4}r_{i}.
\end{equation}
Lemma \ref{lem:orth} gives
\begin{equation}\label{eq:orth2}
\sum_{m,n\in\Fp}\ep(mX+nY+R)=p^{2}\ep(R)\ind_{X=0}\ind_{Y=0}.
\end{equation}
The factors $p^{-2}$ and $p^{2}$ cancel. Hence
\begin{equation}\label{eq:fred}
f_{\chi}(p)=\sum_{r_{1},r_{2},r_{3},r_{4}\in\Fps}
\chibar(r_{1}r_{2})\chi(r_{3}r_{4})\,\ep(r_{1}+r_{2}+r_{3}+r_{4})\,N(r),
\end{equation}
where
\begin{equation}\label{eq:N2}
N(r)=\sum_{\substack{a,b,c,d\in\Fps\\ r_{1}a+r_{2}b+r_{3}c+r_{4}d=0\\
r_{1}a^{2}+r_{2}b^{2}+r_{3}c^{2}+r_{4}d^{2}=0}}
\chi(a)\chi(b)\chibar(c)\chibar(d).
\end{equation}
This proves the reduction \eqref{eq:thm1}--\eqref{eq:N}.

For the exact evaluation we now introduce one more coefficient, rather than apply
unrestricted character orthogonality to the constrained sum \eqref{eq:N2}.

\subsection*{Proof of Theorem 1}

Write $d=\ord(\chi)$ and define
\[
\delta_{j}=\ind_{\chi^{j}=\chi_{0}}\quad (j=2,3,4),\qquad A=|J(\chi,\chi^{2})|^{2}.
\]
Since $d>1$, we have $\delta_{2}=\ind_{d=2}$, $\delta_{3}=\ind_{d=3}$ and
$\delta_{4}=\ind_{d\in\{2,4\}}$.

\smallskip
\noindent\textit{Step 1. Introducing the linear coefficient.}\quad
Set
\[
T(m,n,\ell)=\sum_{a\in\Fps}\chi\bigl(ma^{3}+na^{2}+\ell a\bigr),
\]
and put
\begin{equation}\label{eq:MB}
M=\sum_{m,n,\ell\in\Fp}|T(m,n,\ell)|^{4},\qquad
B=\sum_{m,n\in\Fp}|T(m,n,0)|^{4}.
\end{equation}
For $\ell\not=0$, $T(m,n,\ell)=\chi(\ell)\,C(m/\ell,n/\ell;\chi;p)$. As
$(m,n)\mapsto(m/\ell,n/\ell)$ is a bijection and $|\chi(\ell)|=1$,
\begin{equation}\label{eq:Mrel}
M=(p-1)f_{\chi}(p)+B.
\end{equation}

\smallskip
\noindent\textit{Step 2. The terms with repeated evaluation points.}\quad
For fixed $a,b,c,d\in\Fps$, write
\[
H(a,b,c,d)=\sum_{m,n,\ell\in\Fp}\chi\bigl(F(a)\bigr)\chi\bigl(F(b)\bigr)
\chibar\bigl(F(c)\bigr)\chibar\bigl(F(d)\bigr),
\]
where $F=F_{m,n,\ell}$. Then $M=\sum_{a,b,c,d}H(a,b,c,d)$. If $a=b=c=d$, Lemma
\ref{lem:vander} gives $H=p^{2}\sum_{u\in\Fp}|\chi(u)|^{4}=p^{2}(p-1)$. There are $p-1$
such tuples, so their contribution is $p^{2}(p-1)^{2}$.

If exactly two distinct values occur with multiplicities two and two, there are three
disjoint patterns. For $a=c\not=b=d$ or $a=d\not=b=c$, Lemma \ref{lem:vander} gives
\[
H=p\sum_{u,v\in\Fp}|\chi(u)|^{2}|\chi(v)|^{2}=p(p-1)^{2}.
\]
For $a=b\not=c=d$, it gives
\[
H=p\,\Bigl|\sum_{u\in\Fp}\chi^{2}(u)\Bigr|^{2}=p(p-1)^{2}\delta_{2}.
\]
Each pattern contains $(p-1)(p-2)$ ordered tuples. If two distinct values occur with
multiplicities three and one, or if exactly three distinct values occur, there is a value
occurring only once. Lemma \ref{lem:vander} makes its evaluation variable independent of
the others; summing its non-principal character factor gives zero. Thus all
repeated-point contributions add up to
\begin{equation}\label{eq:Mrep}
M_{\mathrm{rep}}=p^{2}(p-1)^{2}+(2+\delta_{2})\,p(p-1)^{3}(p-2).
\end{equation}

\smallskip
\noindent\textit{Step 3. Four distinct evaluation points.}\quad
Suppose $a,b,c,d$ are pairwise distinct. Apply Lemma \ref{lem:gauss} directly to the four
factors defining $H$. Orthogonality in $m,n,\ell$ yields
\begin{equation}\label{eq:H4}
H(a,b,c,d)=p\sum_{\substack{r_{1},r_{2},r_{3},r_{4}\in\Fps\\
\sum r_{i}t_{i}=\sum r_{i}t_{i}^{2}=\sum r_{i}t_{i}^{3}=0}}
\chibar(r_{1}r_{2})\chi(r_{3}r_{4}),
\qquad (t_{1},t_{2},t_{3},t_{4})=(a,b,c,d).
\end{equation}
The prefactor is $p^{3}/p^{2}=p$. By Lemma \ref{lem:vander}, the solutions are
$r_{i}=\lambda/(t_{i}P'(t_{i}))$ for $\lambda\in\Fps$, where
$P(X)=(X-a)(X-b)(X-c)(X-d)$. In particular,
\[
\frac{r_{3}r_{4}}{r_{1}r_{2}}
=\frac{abP'(a)P'(b)}{cdP'(c)P'(d)}
=\frac{ab(a-b)^{2}}{cd(c-d)^{2}}.
\]
For the last equality, the four cross-difference factors cancel, and both products of
derivatives have the same minus sign. The ratio does not depend on $\lambda$. Define
\[
Q(x,y)=\chi\bigl(xy(x-y)^{2}\bigr),\qquad x,y\in\Fps,\ x\not=y.
\]
We obtain
\begin{equation}\label{eq:HM}
H(a,b,c,d)=p(p-1)\,Q(a,b)\,\overline{Q(c,d)},\qquad
M=M_{\mathrm{rep}}+p(p-1)S,
\end{equation}
where
\begin{equation}\label{eq:Sdef}
S=\sum_{\substack{a,b,c,d\in\Fps\\ a,b,c,d\ \text{pairwise distinct}}}
Q(a,b)\,\overline{Q(c,d)}.
\end{equation}
For $p=3$ this sum is empty; the following identities remain valid.

\smallskip
\noindent\textit{Step 4. Evaluating the distinct-point sum.}\quad
Since $Q(x,y)=Q(y,x)$ and $|Q(x,y)|=1$ for $x\not=y$, put
\[
U_{x}=\sum_{\substack{y\in\Fps\\ y\not=x}}Q(x,y),\qquad V=\sum_{x\in\Fps}U_{x}.
\]
Substituting $y=xt$ gives
\begin{equation}\label{eq:UV}
U_{x}=\chi^{4}(x)\,J(\chi,\chi^{2}),\qquad
\sum_{x}|U_{x}|^{2}=(p-1)A,\qquad
V=(p-1)\delta_{4}J.
\end{equation}
Expand $|V|^{2}$ as a sum over ordered pairs $(a,b)$ and $(c,d)$, each having distinct
nonzero entries. If they have two entries in common, $(c,d)$ is either $(a,b)$ or
$(b,a)$; the contribution is $2(p-1)(p-2)$. If they have exactly one entry in common,
there are four mutually exclusive ways to specify its positions. For instance, the part
with $a=c$ is
\[
\sum_{a}\sum_{\substack{b,d\not=0,\,a\\ b\not=d}}Q(a,b)\overline{Q(a,d)}
=\sum_{a}|U_{a}|^{2}-\sum_{a\not=b}|Q(a,b)|^{2}
=(p-1)A-(p-1)(p-2).
\]
The other three parts have the same value by the symmetry of $Q$. The remaining terms have
four distinct entries and sum to $S$. Consequently,
\begin{equation}\label{eq:S}
\begin{aligned}
S&=|V|^{2}-4\bigl((p-1)A-(p-1)(p-2)\bigr)-2(p-1)(p-2)\\
&=(p-1)\Bigl(\bigl((p-1)\delta_{4}-4\bigr)A+2(p-2)\Bigr).
\end{aligned}
\end{equation}
Combining \eqref{eq:Mrep}, \eqref{eq:HM} and \eqref{eq:S}, we have
\begin{equation}\label{eq:M}
\begin{aligned}
M={}&p^{2}(p-1)^{2}+(2+\delta_{2})\,p(p-1)^{3}(p-2)\\
&+p(p-1)^{2}\Bigl(\bigl((p-1)\delta_{4}-4\bigr)A+2(p-2)\Bigr).
\end{aligned}
\end{equation}

\smallskip
\noindent\textit{Step 5. Removing the zero linear coefficient.}\quad
We have
\[
T(m,n,0)=\sum_{a\in\Fps}\chi^{2}(a)\chi(ma+n).
\]
If $m=0$ and $n\not=0$, this equals $(p-1)\delta_{2}\chi(n)$; if $n=0$ and $m\not=0$, it
equals $(p-1)\delta_{3}\chi(m)$. The term $m=n=0$ is zero. If $mn\not=0$, the
substitution $a=-nt/m$ gives
\[
T(m,n,0)=\chi^{2}(-n/m)\chi(n)J(\chi^{2},\chi).
\]
The Jacobi sum is symmetric in its two characters by $t\mapsto1-t$, so the fourth absolute
power is $A^{2}$. Counting the parameter pairs in the three nonzero cases gives
\begin{equation}\label{eq:B}
B=(\delta_{2}+\delta_{3})(p-1)^{5}+(p-1)^{2}A^{2}.
\end{equation}

\smallskip
\noindent\textit{Step 6. Substitution for each character order.}\quad
Equations \eqref{eq:Mrel}, \eqref{eq:M} and \eqref{eq:B} give the uniform identity
\begin{equation}\label{eq:f}
\begin{aligned}
f_{\chi}(p)={}&p^{2}(p-1)+(2+\delta_{2})\,p(p-1)^{2}(p-2)\\
&+p(p-1)\Bigl(\bigl((p-1)\delta_{4}-4\bigr)A+2(p-2)\Bigr)
-(\delta_{2}+\delta_{3})(p-1)^{4}-(p-1)A^{2}.
\end{aligned}
\end{equation}
For $d=2$, the values are $(\delta_{2},\delta_{3},\delta_{4},A)=(1,0,1,1)$, and \eqref{eq:f}
becomes
\[
f_{\chi}(p)=(p-1)\Bigl[p^{2}+3p(p-1)(p-2)+p(3p-9)-(p-1)^{3}-1\Bigr]
=2p(p-1)(p^{2}-p-3).
\]
For $d=3$, they are $(0,1,0,1)$, so
\[
f_{\chi}(p)=(p-1)\Bigl[p^{2}+2p(p-1)(p-2)+p(2p-8)-(p-1)^{3}-1\Bigr]
=p(p-1)(p^{2}-7).
\]
For $d=4$, they are $(0,0,1,p)$, so
\[
f_{\chi}(p)=p(p-1)\Bigl[2(p-1)(p-2)+p(p-5)+2(p-2)\Bigr]
=3p^{2}(p-1)(p-3).
\]
Finally, if $d\notin\{2,3,4\}$, they are $(0,0,0,p)$, and
\[
f_{\chi}(p)=p(p-1)\Bigl[2(p-1)(p-2)-4p+2(p-2)\Bigr]=2p^{2}(p-1)(p-4).
\]
These are precisely the four branches of \eqref{eq:main}. This proves Theorem 1.

\subsection*{Proof of Corollary 1}

Under the additional assumption $d\notin\{2,3,4\}$, the fourth branch of \eqref{eq:main} gives
\begin{equation}\label{eq:asym}
f_{\chi}(p)=2p^{2}(p-1)(p-4)=2p^{4}-10p^{3}+8p^{2}=2p^{4}-10p^{3}+O(p^{2}).
\end{equation}
The implied constant can be taken to be $8$, independently of $\chi$. Dividing by $p^{2}$
yields the normalized assertion. For the exceptional orders, the other branches give
\begin{equation}\label{eq:branches}
\begin{aligned}
d=2:\quad & f_{\chi}(p)=2p^{4}-4p^{3}-4p^{2}+6p,\\
d=3:\quad & f_{\chi}(p)=p^{4}-p^{3}-7p^{2}+7p,\\
d=4:\quad & f_{\chi}(p)=3p^{4}-12p^{3}+9p^{2}.
\end{aligned}
\end{equation}
The leading coefficients are therefore $2$, $1$, $3$, respectively.

\begin{remark}[Structure of the exceptional cases]
Formula \eqref{eq:f} identifies all sources of exceptional behaviour. The repeated-point pattern
$a=b\not=c=d$ introduces $\delta_{2}$. The sum $V$ introduces $\delta_{4}$, since its
scaling factor is $\chi^{4}(x)$. The zero-linear-coefficient terms introduce $\delta_{2}$
and $\delta_{3}$. Finally, $A=|J(\chi,\chi^{2})|^{2}$ is $1$ for orders $2$ and $3$, and
$p$ otherwise. Thus only orders $2,3,4$ differ from the generic case. This conclusion
follows from the explicit evaluations above, including the nonzero contribution of four
distinct points, rather than from an assumption that all non-diagonal blocks vanish.
\end{remark}

%=====================================================================
\section{Numerical verification}
%=====================================================================

Let $g$ be a primitive root mod $p$. Every Dirichlet character mod $p$ is of the form
\[
\chi_{k}\bigl(g^{t}\bigr)=\exp\Bigl(\frac{2\pi ikt}{p-1}\Bigr),\qquad 0\le k\le p-2,
\]
and its order is $d=(p-1)/\gcd(k,p-1)$. The character $\chi_{0}$ is principal and is
excluded from Theorem 1. Direct evaluation of \eqref{eq:def} requires $O(p^{3})$ operations per
character, and $O(p^{4})$ operations per prime when all characters are evaluated
separately.

We independently recomputed \eqref{eq:def} for all $24$ odd primes $p\le97$ and all $1010$
non-principal characters. For each $(m,n)$, the nonzero values of $ma^{3}+na^{2}+a$ were
counted according to their discrete logarithms; a discrete Fourier transform then
evaluated all characters simultaneously. Complex conjugation changes the sign of the
Fourier exponent but leaves the fourth absolute moment unchanged. The largest absolute
discrepancy between these floating-point values and \eqref{eq:main} was less than $3\times10^{-7}$.
These computations provide numerical checks; the proof of \eqref{eq:main} is the exact argument in
Section 3, not rounding of numerical output.

Table 1 lists the results for $3\le p\le47$. For a fixed $(p,d)$, all characters of order
$d$ give the same value, as proved in Theorem 1. The table covers $14$ primes and $298$
non-principal characters. The remaining primes checked were $53$, $59$, $61$, $67$, $71$,
$73$, $79$, $83$, $89$, $97$.

\begin{table}[htbp]
\centering
\caption{Fourth power means for $3\le p\le47$: 14 primes and 298 non-principal characters.}
\small
\begin{tabular}{r@{\hspace{2.6em}}l@{\hspace{2.6em}}l}
\toprule
$p$ & $d=\ord(\chi)$ & $f_{\chi}(p)$\\
\midrule
3  & 2 & 36\\
5  & 2 / 4 & 680 / 600\\
7  & 2 / 3 / 6 & 3\,276 / 1\,764 / 1\,764\\
11 & 2 / 5, 10 & 23\,540 / 16\,940\\
13 & 2 / 3 / 4 / 6, 12 & 47\,736 / 25\,272 / 60\,840 / 36\,504\\
17 & 2 / 4 / 8, 16 & 146\,336 / 194\,208 / 120\,224\\
19 & 2 / 3 / 6, 9, 18 & 231\,876 / 121\,068 / 194\,940\\
23 & 2 / 11, 22 & 509\,036 / 442\,244\\
29 & 2 / 4 / 7, 14, 28 & 1\,313\,816 / 1\,836\,744 / 1\,177\,400\\
31 & 2 / 3 / 5, 6, 10, 15, 30 & 1\,724\,220 / 887\,220 / 1\,556\,820\\
37 & 2 / 3 / 4 / 6, 9, 12, 18, 36 & 3\,540\,456 / 1\,814\,184 / 5\,026\,968 / 3\,252\,744\\
41 & 2 / 4 / 5, 8, 10, 20, 40 & 5\,369\,360 / 7\,665\,360 / 4\,975\,760\\
43 & 2 / 3 / 6, 7, 14, 21, 42 & 6\,512\,436 / 3\,326\,652 / 6\,057\,324\\
47 & 2 / 23, 46 & 9\,335\,516 / 8\,738\,804\\
\bottomrule
\end{tabular}

\smallskip
{\footnotesize Entries in the column $d=\ord(\chi)$ separated by a slash correspond to the
successive values of $f_{\chi}(p)$; entries separated by a comma share the same value.}
\end{table}

Finally, the principal character is not covered by \eqref{eq:main}, but its fourth power mean also has
a polynomial expression. Write
\[
R(m,n)=\#\bigl\{a\in\Fps:\ ma^{2}+na+1=0\bigr\}.
\]
Then $C(m,n;\chi_{0};p)=p-1-R(m,n)$. Let $N_{j}$ denote the number of pairs $(m,n)$ for
which $R(m,n)=j$. Two distinct nonzero roots $u,v$ determine the unique polynomial
$(1-X/u)(1-X/v)$, so $N_{2}=(p-1)(p-2)/2$. Exactly one nonzero root arises either from a
nonconstant linear polynomial, giving $p-1$ possibilities, or from a quadratic polynomial
with a double root, giving another $p-1$. Hence
\[
N_{0}=\frac{p^{2}-p+2}{2},\qquad N_{1}=2(p-1),\qquad N_{2}=\frac{(p-1)(p-2)}{2}.
\]
It follows that
\begin{equation}\label{eq:principal}
\begin{aligned}
f_{\chi_{0}}(p)&=N_{0}(p-1)^{4}+N_{1}(p-2)^{4}+N_{2}(p-3)^{4}\\
&=p^{6}-8p^{5}+34p^{4}-92p^{3}+153p^{2}-138p+50.
\end{aligned}
\end{equation}
For $p=3,5,7,11,13$, this gives $68$, $3\,560$, $39\,852$, $875\,540$ and $2\,649\,528$,
respectively, agreeing with the values obtained directly.

%=====================================================================
\section*{Competing interests}
%=====================================================================

The authors declare that they have no competing interests.

\section*{Acknowledgments}

The authors would like to thank the referees for their very helpful and detailed comments,
which have significantly improved the presentation of this paper.

%=====================================================================

\end{document}